\documentclass[notitlepage, 11pt]{article}
\usepackage{hyperref, amsmath, amsfonts, enumitem, amsthm, xcolor, parskip}
\usepackage{amssymb} 
\usepackage[nottoc]{tocbibind}
\usepackage{tikz-cd}
\usepackage{subfiles}

\hypersetup{             
    colorlinks=true,                
    breaklinks=true,                
    urlcolor= black,                
    linkcolor= blue,                     
    bookmarksopen=false,
    filecolor=black,
    citecolor=blue,
    linkbordercolor=blue
}
\setlist[itemize]{itemsep=0pt, parsep=0pt}
\setlist[enumerate]{itemsep=0pt, parsep=0pt}
\theoremstyle{definition}
\newtheorem{defn}{Definition}
\newtheorem{eg}{Example}[section]
\theoremstyle{plain}
\newtheorem{thm}[defn]{Theorem}

\newtheorem{prop}[defn]{Proposition}
\newtheorem{cor}{Corollary}[defn]
\theoremstyle{remark}

\newcommand{\R}{\mathbb{R}}

\newcommand{\C}{\mathbb{C}}
\newcommand{\bH}{\mathbb{H}}
\newcommand{\N}{\mathbb{N}}
\newcommand{\Z}{\mathbb{Z}}
\newcommand{\F}{\mathbb{F}}

\title{Fields With Finitely Many Non-Commutative Division Algebras Over Them.}
\author{Snehinh Sen\footnote{\textsc{Theoretical Statistics and Mathematics Unit\\ Indian Statistical Institute, Kolkata Centre, \\203, B.T. Road, Kolkata - 700107}\\
E-mail: snehinh1math@gmail.com}}
\setlist[enumerate, 1]{label={(\roman*)}}

\begin{document}

\maketitle

\begin{abstract}
We classify fields having finitely many finite non-commutative (not necessarily central) division algebras over them. In the process, we introduce the notion of \textit{anti-closure} of a field and also make comments on fields having a linear lattice of finite field extensions over them. 
\end{abstract}



\section{Introduction}\label{s1}

Frobenius showed in \cite{fro78} that the only finite division algebras over $\R$, upto isomorphism, are $\R$, $\C$ and $\bH$, the Hamiltonian Quaternions. In fact, this result generalizes to arbitrary real closed fields (see, for eg., \cite{ger60}). As a consequence, if a field has finitely many isomorphism classes of finite (possibly commutative) division algebras over it, it is either real closed or algebraically closed. However, this categorization excludes fields such as the finite fields, over which every finite division algebra is commutative.

The main goal of this article is to classify fields which have only finitely many finite non-commutative division algebras upto isomorphism over them. It turns out, interestingly, that beyond the aforementioned class, real closed fields are the only other possibility yet again. Here is the precise statement. 

\textbf{Theorem.}\emph{\label{t02} 
A field $K$ having only finitely and positively many nonisomorphic noncommutative finite division algebras over it is real closed. }

To the best of the knowledge of the author, this result is not present in the literature. The article is structured as follows. Section \ref{s2} is a recap of results from the theory of Brauer groups and division algebras. Section \ref{s3} gives a short proof of the generalized Frobenius Theorem. Sections \ref{s4} and \ref{s5} deal with the notion of \textit{anti-closure}, its properties and its connection to fields having a linear lattice of finite extensions over them. Sections \ref{s6} and \ref{s7} is devoted to preparing machinery and settling the claim made above.

\section{Division Algebras and Brauer Group} \label{s2}

We briefly recall some basic properties of division algebras and Brauer groups without proof. Suitable references are \cite{gui18}, \cite{ser79}, \cite{ser97} and \cite{stack23}.

\begin{defn}\label{d:Divalg}
Let $K$ be a field and $D$ be a finite division algebra over $K$.
\begin{enumerate}
\item $r_K(D)$ denotes the \textit{rank} of $D$ over $K$, that is $dim_K(D)$.
\item $Z(D)$ denotes the center of $D$.
\item It is known (see, for eg., \cite{stack23}) that $r_{Z(D)}(D)$ is a square of an integer. Its positive square root, denoted by $\deg(D)$, is called the \textit{degree} of $D$.
\end{enumerate}
\end{defn}

\begin{defn}\label{d:Brauer} 
Let $K$ be a field and $L/K$ be a Galois field extension.  
\begin{enumerate}
\item $Br(L/K)$ denotes the \textit{Relative Brauer Group} or, equivalently, the second cohomology group $H^2(G(L/K), L^*)$.
\item $Br(K) = H^2(/K)$ denotes the \textit{Brauer Group} of the field $K$. 
\item $S(K)$ denotes the classes of all finite non-commutative division algebras over $K$ upto isomorphism.
\end{enumerate}
\end{defn}

\begin{prop}\label{p1} 
Let $K$ be a field, $L/K$ a finite Galois extension, $D$ a finite division algebra over $K$ and $\overline{K}$ an algebraic closure or $K$ containing $L$.
\begin{enumerate}
\item We have the inflation-restriction exact sequence $$0 \to Br(L/K) \to Br(M/K) \to Br(M/L)$$ where $M/K$ is any Galois extension containing $L$. 
\item We have an exact sequence $$0 \to Br(L/K) \to Br(K) \to Br(L)$$ where the first map is essentially inclusion and the second map is given by $A\mapsto A\otimes_K L$ upto Morita equivalence. 
\item Suppose $[D]\in Br(K)$ maps to $0$ in the above sequence (in other words, $L$ is a \emph{splitting field }of $D$), then $\deg(D)|[L:K]$.
\item If $g = |G(L/K)|$, then $g.Br(L/K)=0$.
\item $Br(K)$ is a torsion group which is a union of all $Br(L/K)$ as $L/K$ varies over all finite Galois extensions of $K$.
\item As sets, there is a natural bijection between $Br(K)$ and the set of isomorphism classes of finite $K-$central division algebras $D$.
\item As sets there is a natural bijection between $S(K)$ and the disjoint union 
\[ \coprod_{L} \left(Br(L)\setminus \{[0]\}\right)\]
\end{enumerate}
where $L$ runs over finite extensions of $K$ inside $\bar{K}$ which are not conjugates. 
\end{prop}

\begin{cor}\label{p1.1}
Let $K$ be a field. Then $K$ has only finitely many finite noncommutative division algebras upto isomorphism over it if and only if 
\begin{enumerate}
\item For each $L/K$ finite, $Br(L)$ is finite, and
\item For all but finitely many such extensions, $Br(L)=0$. 
\end{enumerate}
\end{cor}

\begin{proof}
This is immediate from parts (vi) and (vii) of Proposition \ref{p1}. 
\end{proof}

\begin{cor}\label{p3.1}
Let $K$ be a separably closed field. Then $Br(K)=0$.
\end{cor}

\begin{proof}
As $K^{sep}=K$, $Br(K)=Br(K/K) = 0$.
\end{proof}

\begin{prop}\label{p2} 
Let $K$ be a field and $L/K$ be a cyclic extension. Then $Br(L/K) = K^*/N(L^*)$ where $N = N_{L/K}$ is the norm map.
\end{prop}

\begin{proof}
Follows from \cite{ser79}, Chapter VIII, Proposition 6. 
\end{proof}

\begin{cor}\label{cp2.1}
Let $K$ be a perfect field with characteristic $p$ and let $L/K$ be a degree $p$ extension. Then the norm map $N_{L/K}$ is surjective. Consequently, if $M/K$ is cyclic and has degree a power of $p$, then $Br(M/K)=0$. 
\end{cor}

\begin{proof}
The first part follows as $K^* = (K^*)^p \subseteq N_{L/K}(L^*)$. For the Galois case, let $K=K_0\subseteq K_1\ldots K_n = M$ denote a tower of extensions such that each $[K_i:K_{i-1}] = p$. Then each $N_{K_i/K_{i-1}}$ is surjective, proving the surjectivity of $N_{M/K}$. Rest follows from Proposition \ref{p2}.  
\end{proof}

\begin{cor}\label{cp2.2}
Let $L/K$ be a cyclic extension such that $N_{L/K}$ is surjective and $Br(L)=0$. Then $Br(K)=0$.
\end{cor}

\begin{proof}
The exact sequence of Propsotion \ref{p1}, part (ii) and Proposition \ref{p2}, yield that $0\to Br(K)\to 0$ is exact, showing that $Br(K)=0$.
\end{proof}

\section{Frobenius Theorem for Real Closed Fields} \label{s3}

We shall now give a short proof of Frobenius Theorem for real closed fields and provide motivation for our main problem through a corollary. This result is already known (see, for instance, \cite{ger60}) but is proved nevertheless as its reinterpretation plays a key role for us. For psychological reasons, let $R$ be a real closed field and $C$ be an algebraic closure of $R$. It is known that $C = R(i)$ where $i^2=-1$. Likewise, let $H$ denote the rank $4$ $R$-central division algebra defined as follows.

\begin{enumerate}
\item $H$ has an $R-$ basis denoted by $1,i,j,k$ and $Z(H)=R1$.
\item $i^2=j^2=k^2=ijk = - ikj = -1$.
\end{enumerate}

Similar to the real case, one may verify that $H$ is a division algebra. In fact, the inverse of a non-zero $\alpha = a+bi+cj+dk$ is given by $\frac{1}{N(\alpha)}(a-bi-cj-dk)$ where $N(\alpha) = a^2+b^2+c^2+d^2>0$ as $R$ is real closed. 

\begin{thm}[Frobenius Theorem for Real Closed Fields]\label{t1} 
Suppose $R$ is a real closed field. Then the only finite division algebras over $R$, upto isomorphism, are $R, C$ and $H$ as constructed above. 
\end{thm}

\begin{proof}
Note that $Br(R) = Br(C/R) = R^*/N(C^*) \cong R^*/(R^{*})^2\cong C_2$. So the only two $R-$central division algebras upto isomorphism are $R$ and $H$. However, any other finite $R-$division algebra must be $C-$central. By Corollary \ref{p3.1}, the only possibility here is $C$ itself. This completes the proof.
\end{proof}

\begin{cor}\label{ct1.1}
A field $K$ has only finitely many non-isomorphic finite division algebras over it if and only if $K$ is real closed or algebraically closed. 
\end{cor}

\begin{proof}
Firstly, as every finite field extension of $K$ is a division algebra over it, $K$ should have only finitely many finite non-conjugate field extensions. By Artin-Schreier Theorem, this is possible only if $K$ is real closed or algebraically closed. Converse is Theorem \ref{t1} and Corollary \ref{p3.1}.
\end{proof}

This corollary is the prime motivation for our considerations - what happens if we restrict our attention only to non-commutative division algebras? 

\section{Anti-Closure and M-Groups} \label{s4}

The notion of anti-closure of a field, though not standard, will pop up in our analysis. Firstly, we define what we mean by the anti-closure.

\begin{defn}\label{d:anticlosure}
Let $K$ be a field and $\overline{K}$ be an algebraic closure of $K$. The \textit{anti-closure} of $K$ in $\overline{K}$ is defined as 
\[ K' := \bigcap_{L} L\]where $L$ runs over all finite non-trivial extensions of $K$ in $\bar{K}$.
\end{defn}

The reason for this nomenclature is that the algebraic closure satisfies $\overline{K} = \bigcup_L L$ where the union is over a similar collection. If $K$ is algebraically closed, we define $K'$ to be $K$. It must be noted that $K'/K$ is a non-trivial field extension if and only if $K$ has a unique minimal extension over it. We now state a few properties of fields for which $K'\neq K$. Henceforth, we fix the algebraic closure and consider finite subextensions only.

\begin{prop}\label{p4:ac}
    Let $K$ be a field such that $K'\neq K$. Then we have the following properties. 
    \begin{enumerate}
     \item If $L/K$ is a non-trivial extension of $K$, then $K'\subseteq L$.
    \item $K'/K$ is finite and normal.
    \item $K$ is either perfect or separably (but not algebraically) closed.
    \item $[K':K]$ is a prime number $p$.
    \item If $K$ is not perfect, then $char(K)=[K':K]$
    \item Let $p$ be the prime above. Then every extension of $K$ is a $p-$extension.
\end{enumerate}
\end{prop}

\begin{proof}
The definition of $K'$ implies (i) and that $K'/K$ is a finite extension contained in all of its Galois conjugates. So $K'/K$ is normal, proving (ii).

 For (iii), suppose $K$ is not perfect. Then $K$ has a non-trivial purely inseparable extension, say $L/K$. As $K'\subseteq L$, we get that $K'/K$ is purely inseparable. Hence, every non-trivial extension of $K$, by the virtue of it containing $K'$, is inseparable, showing that $K$ is separably closed.  

For (iv) and (v), if $K$ is perfect, then $K'/K$ is Galois and $G(K'/K)$ is a finite group with no proper subgroups. Hence, $G(K'/K)\cong C_q$ for some prime $q$, proving our claim. If $K$ is not perfect, then let $\alpha\in K\setminus K^p$, where $p=char(K)$. Then $K(\alpha^{1/p})$, being a minimal extension of $K$, is equal to $K'$. Thus $[K':K]=p$. This immediately implies (iv) and (v). 

Finally, if $K$ is separably closed, (vi) is immediate. Otherwise, $K$ is perfect by (iv). Let $L/K$ be a finite non-trivial extension of $K$ and $N/K$ be its normal closure in $\overline{K}$. Let $p=[K':K]$ and $G=G(N/K)$. Let $G_p$ be a Sylow p-subgroup of $G$ and $M$ its fixed field. Then $([M:K],p)=1$. Hence $M$ does not contain $K'$, implying that $M=K$ and $L/K$ is a $p$-extension. 
\end{proof}

Perfect fields with $K'\neq K$ are very rare. Here are a few examples.

\begin{eg}\label{e1}
\begin{enumerate}
\item Let $K$ be a real closed field. Then $K'$ is simply an algebraic closure of $K$. Conversely, if $K'$ is algebraically closed, then $K$ is real closed or algebraically closed (by Artin-Schreier Theorem).
\item Let $K$ be a perfect field with absolute Galois group isomorphic to $\Z_p$ for some $p$ prime. For example, let $K$ denote the fixed field of the $p-$part $\Z_p\leq \hat{\Z}$ for the Galois extension $\overline{\F}_q/\F_q$ (which has Galois group $\hat{\Z}$) where $q$ is any prime power, Likewise, let $K=\C((T^{1/n} : (p,n)=1))$. Then $K'$ is the fixed field of $p\Z_p$ and is not equal to $K$.
\end{enumerate}
\end{eg}

As we shall see soon, these are in fact the only examples. Now fields with $K'\neq K$ are very much similar to real closed fields as we shall show next. 

\begin{prop}\label{p5:acrc}
    Let $K$ be a field such that $K'\neq K$. Let $p=[K':K]$
    \begin{enumerate}
        \item The degrees of irreducible polynomials in $K[X]$ are powers of $p$.
        \item Every polynomial in $K[X]$ having degree coprime to $p$ has a root in $K$.
        \item $K$ contains all the $p-$th roots of unity.
        \item Every element in $K$ is a $p^{th}$ power of an element in $K'$ and $K'$ contains all the $p^{th}$ roots of elements of $K$.
        \item If $char(K)\neq p$ or $K$ is separably closed, then there is an element $\alpha \in K'$ such that $\alpha^p \in K$ and $K'=K(\alpha)$. So $K$ is not closed under taking $p^{th}$ roots in this case. 
    \end{enumerate}
\end{prop}

\begin{proof}
Let $f(X)$ be an irreducible polynomial and $L/K$ be a splitting field. Then $\deg(f)$ divides the order of $[L:K]$, which is a power of $p$ by Proposition \ref{p4:ac}. Hence $\deg(f)$ is a power of $p$, proving (i). Now suppose $f(X)$ has degree coprime to $p$. Then $f$ has an irreducible factor with degree coprime to $p$. (i) would imply that such a factor is linear, proving (ii). Applying these on $X^{p-1}+\ldots+1$ and $X^p-a$ yield (iii) and (iv) immediately. (v) is a consequence of Kummer Theory and basic field theory, along with part (iv).
\end{proof}

In analogy to the notion of anti-closure, we define a concept in group theory which shall help us study such field extensions better. Such groups, though existent in literature, have usually been studied anonymously. We give them a name for an easier reference.

\begin{defn}\label{d:Mgroups}
An \textit{M-group} is a Hausdorff topological group $G$ with a unique closed maximal subgroup.
\end{defn}

Note that if $K$ is perfect and $G=G(/K)$, then $K'\neq K$ if and only if $G$ is an $M-$group. Furthermore, we have the following proposition.

\begin{prop}\label{p6} 
Every $M$-group $G$ is monothetic. 
\end{prop} 

\begin{proof}
Let $x$ be an element of $G$ which is not in its maximal closed subgroup. Then $\overline{\langle x\rangle}$ is a closed subgroup of $G$ which is not contained in the unique maximal subgroup. Hence, it must be the whole group $G$.
\end{proof}

\begin{cor}\label{cp6.1}
Every finite (resp. profinite) $M$-group is cyclic (resp. procyclic). 
\end{cor}

The following result satisfactorily classifies the absolute Galois group of all perfect fields such that $K'\neq K$. 

\begin{thm}\label{t3} 
Let $K$ be perfect with $K'\neq K$ and $G=G(/K)$. Then $G$ is either isomorphic to $C_2$ or the additive group $\Z_p$ where $p=[K':K]$.
\end{thm}

\begin{proof}
Note that $G$ is finite if and only if $K$ is real closed and $G\cong C_2$. So let $G$ be infinite. By Proposition \ref{p4:ac} and Corollary \ref{cp6.1}, $G$ is a pro-$p$ procyclic group. Let $S = \{\log_p([G:N]) : N\text{ is a normal, finite index subgroup of }G\}$. As $G$ is infinite, $S$ is cofinal in $\N$ and hence $G \cong  \hat{G} = \varprojlim G/N = \Z_p$. 
\end{proof}

Before moving forward, we comment on the case when $K$ is separably closed.

\begin{prop}\label{pinsep} 
Let $K$ be a separably closed field of characteristic $p\neq 0$. Then the following are equivalent.
\begin{enumerate}
\item $K'\neq K$.
\item $[K:K^p]=p$.
\item $K/K^p$ is simple and nontrivial.
\end{enumerate}
\end{prop}

\begin{proof}
Assume (i). Then $K'=K(a^{1/p})$ for some $a\in K\setminus K^p$. Clearly, for each $b\in K$, $b^{1/p} \in K'=K(a^{1/p})$ implies that $b\in K^p(a)$. So $K=K^p(a)$. proving (iii). Assuming (iii), if $K=K^p(b)$, then the minimal polynomial of $b$ in $K$ divides $g(X)=X^p-b^p = (X-b)^p$. As $b\notin  K^p$, $k=p$ is the smallest positive integer such that $b^k\in K^p$. Thus, $g(X)$ is the minimal polynomial of $b$ over $K^p$, proving (ii). Assuming (ii), if $a,b\in K\setminus K^p$, then $K=K^p(a)=K^p(b)$. So $K(a^{1/p}) = K(b^{1/p})$, proving that $K$ has a unique extension of degree $p$. As $K$ is separably closed, every minimal extension is of degree $p$. Hence, $K$ has a unique minimal extension, proving (i).
\end{proof}

\section{Fields with Linear Lattices} \label{s5}

In this section, we shall classify the structure and behaviour of fields satisfying a peculiar condition with respect to their lattice of extensions.

\begin{defn}\label{d:ll}
Let $K$ be a field. We say that $K$ has a \textit{linear lattice over it} or satisfies the linear lattice condition if there are field extension $K^{(n)}$, $n\geq 0$ (possibly up to a finite index) over $K$ such that
\begin{enumerate}
\item $K= K^{(0)}$.
\item $K^{(i)}\subseteq K^{(i+1)}$ for each $i\geq 0$.
\item For every $L/K$ finite, there is an $i\geq 0$ such that $L=K^{(i)}$.
\end{enumerate}
We would say that this is a \textit{linear lattice of distinct field extensions} over $K$, if $K^{(i)}\neq K^{(i+1)}$ for each $i$. 
\end{defn}

An independent theory of such fields can be developed. However, we would use the results of the previous section to classify such fields.

\begin{thm}\label{t4} 
Let $K$ be a perfect field which is not algebraically closed. Then the following are equivalent.
\begin{enumerate}
\item $K$ has a linear lattice over it.
\item $K'\neq K$.
\item $G(/K)$ is an $M$-group.
\item $G(/K)$ is isomorphic to either $C_2$ or $\Z_p$ for some prime $p$.
\end{enumerate}
\end{thm}

\begin{proof}
$\text{(ii)}\iff \text{(iii)}\iff \text{(iv)}$ has already been established. Suppose $\text{(i)}$ is assumed. Then, as $K$ is not algebraically closed, there must be a least $i\geq 1$ such that $K^{(i)}\neq K$. We realise that $K'=K^{(i)}\neq K$, proving $\text{(ii)}$. Now, as the closed subgroups of $C_2$ and $\Z_p$ form a linear lattice, we immediately get $\text{(i)}$ by assuming $\text{(iv)}$ via the Galois correspondence.
\end{proof}

We are prepared to make structural comments on fields with such lattices. The analysis would be split into two cases. Let us say that $[K':K]=p>1$.

\subsection{Characteristic of \texorpdfstring{$K$}{K} is \texorpdfstring{$p$}{p}}

We shall use Artin-Schreier Theory to derive the following claim.

\begin{thm}\label{t5}
    Let $K$ be a perfect field such that $[K':K]=char(K)=p$. Let $\mathcal{P}(T)=T^p-T$. There is an element $\alpha_0 \in K$ such that if we inductively define $\alpha_{i+1}$ to be a root of $\mathcal{P}(T)-\alpha_i$ for each $i$, then $K^{(i)}=K(\alpha_i)$ gives the linear lattice of distinct field extensions over $K$. In fact, $\alpha_0$ can be chosen to be any element of $K$ such that $\mathcal{P}(T)-\alpha_0$ is irreducible.
\end{thm}

\begin{proof}
    Existence of such an $\alpha_0$ is asserted by Artin-Schreier Theory. Suppose, for the sake of a contradiction, $i\geq 1$ is smallest such that $\beta = \alpha_{i+1} \in K(\alpha_i)=L$. Then let $\gamma =\alpha_i$. Now let $x \in L$. Let $F = K(\alpha_{i-1})$. By Artin-Schreier Theory, $i>1$. So $x = a_0+a_1\gamma+\ldots+a_{p-1}\gamma^{i-1}$ for some $a_i \in F$. Now, as $K$ is perfect, there are elements $b_i \in F$ such that $a_i = b_i^p$. Thus $y = b_0+b_1\beta+\ldots+b_{p-1}\beta^{i-1}$ satisfies $\mathcal{P}(T)-x$. Hence, by Artin-Schreier Theory, $L$ has no degree $p$ extension. As $G(/L)$ is pro-$p$, we get $L$ is algebraically closed, contradicting the Artin-Schreier Theorem. 
\end{proof}

\textbf{Remark.} As $K$ is perfect with a positive characteristic, the norm map of any extension of $K$ is surjective. By cyclicity, $Br(K^{(i)}/K)=0$ for each $i\geq 0$ by Corollary \ref{cp2.1}. Taking limit on $i$, we obtain that $Br(K)=0$.

\subsection{Characteristic of \texorpdfstring{$K$}{K} is not \texorpdfstring{$p$}{p}}

Here Kummer Theory would be applicable instead of Artin-Schreier Theory. First, we will deal with the norm. 

\begin{prop}\label{p7}
    Let $K$ be a field such that $char(K)$ is not equal to $p=[K':K]$. Then the norm map $N_{K'/K}$ is surjective with the only exception being the case when $K$ is real closed. In that case, $\mathrm{Coker}(N_{K'/K})\cong C_2$.
\end{prop}

\begin{proof}
    Such a $K$ is perfect by Propositon \ref{p4:ac}. For each $a\in K$, either $X^p - a$ is irreducible or it has a root in $K$. In either case, the root lies in $K'$ by Proposition \ref{p5:acrc}. In the first case, the minimal and characteristic polynomial is the same as $X^p-a$. So the norm is $(-1)^{p-1}a$. In the latter case, the characteristic polynomial is $(X-b)^p$, where $b^p = a$. So the norm is $a$. If $p\neq 2$, then this would immediately imply that the norm is surjective. 
    
    Now let $p=2$. Firstly, suppose $-1$ is not a square in $K$. Here $K'=K(i)$ where $i^2 = -1$. For each $a \in K$, $\sqrt{a}=c+di$ for some $c,d\in K$. Thus $a = (c+di)^2$, which would imply that $a=c^2$ or $a=-d^2$. Hence, every element of $K$ is either a square or the negative of a square. If $(-1)$ is in the image of the norm, then we would get that the norm map is surjective. 
    
    Otherwise, this shows that $N(L^*)=(K^*)^2$. So $K^*/N(L^*) \cong C_2$. If $x=s^2, y=t^2$, where $s,t\in K$, then $s+it \in K(i)$ will have norm $s^2+t^2$. So $K$ is pythagorean. For $c+di\in K(i)$, $c,d\in K$, if we set $\alpha = \sqrt{\frac{c+\sqrt{c^2+d^2}}{2}} + i\sqrt{\frac{-c+\sqrt{c^2+d^2}}{2}} \in K(i)$, we get $\alpha^2 = c+di$. Hence, $K'$ has no quadratic extensions, proving that $K'$ is algebraically closed and thus $K$ is real closed.
    
    Finally, suppose $(-1)$ is a square in $K$. Then, $X^2-a$ is irreducible if and only if $X^2+a$ is irreducible. Thus, every element in $K$ would appear as a norm from the first paragraph.
\end{proof}

So, once again, we get that $Br(K) = 0$ with the only exception being the case when $K$ is real closed. We now comment on the lattice above $K$.

\begin{thm}\label{t6}
    Let $K$ be a perfect field with $char K \neq p = [K':K]$. Suppose $(-1)$ is a $p^{th}$ power in $K$. There is an element $\alpha \in K=K^{(0)}$ such that $K^{(i)} = K(\alpha^{p^{-i}})$ gives the linear lattice of distinct field extensions over $K$. In fact, $\alpha$ can be chosen to be any element which is not a $p^{th}$ power in $K$.
\end{thm}

\begin{proof}
    Such an $\alpha$ exists by Kummer Theory. Suppose, on the contrary, $i$ is the least index for which $\alpha$ chosen as above satisfies $\alpha^{p^{-i}} \in K^{(i-1)} = L $. Clearly $i> 1$. Set $F = K^{(i-2)}$, $\beta = \alpha^{p^{-(i-1)}}$ and $\gamma = \alpha^{p^{-i}}$. Then $L = F(\beta)= F(\gamma)$. Now let $N$ be the norm map from $L$ to $F$. Then $N(\gamma) = a$ implies that $N(\beta) = a^p$. But the characteristic polynomial of $\beta$ is of the form $X^p - c$, where $c$ is not a $p^{th}$ power in $F$. But the norm is $a^p = (-1)^{p-1}c$, showing that $c$ is a $p^{th}$ power. This is a contradiction.
\end{proof}

This settles the case when $p$ is odd or when $p=2$ and $-1$ is a square in $K$. If $-1$ is not a square in $K$, there are two possibilities, as can be seen in the proof of Proposition \ref{p7}. Namely, $K$ is either real closed or $(-1)$ is in the image of the norm map $N_{K'/K}$. The first case has a simple lattice structure, with $K^{(1)}=K'$ algebraically closed. Conversely, according to Satz 5 of \cite{art27}, every algebraically closed field of characteristic $0$ contains a real closed subfield. For the second case, we look at the following theorem. 

\begin{thm}\label{t7} 
    Let $K$ be a perfect field with $[K':K]=2$. Suppose $K$ is not real closed and $(-1)$ is not a square in $K$. Then every element of $K$ is a fourth power in $K'$. The linear lattice of distinct field extensions over $K$ is given by $K^{(0)}=K$ and $K^{(n+1)} = K(\alpha_0^{2^{-n}})$ for $n\geq 0$, where $\alpha_0$ is any element of $K'$ such that $N_{K'/K}(\alpha_0)$ is not a square in $K$.
\end{thm}

\begin{proof}
The proof of Proposition \ref{p7} shows that every element of $K$ is either a square or its negative is a square in $K$. Let $i\in K'$ satisfy $i^2=-1$. Now $\beta = u(1+i)$, where $u = \sqrt{\frac{1}{2}}$, satisfies $\beta^2 = i$. As $\sqrt{a}\in K'$ for each $a\in K$, we get that every element of $K$ is a fourth power in $K'$.

For the claim on lattices, we let $L=K'$. We realise that $L$ also has a linear lattice over it and also $(-1)$ is a square in $L$. So, by Theorem \ref{t6}, the linear lattice of distinct field extensions can be described by $L^{(i)} = L(\alpha^{2^{-i}})$ where $\alpha$ is an element of $L$ which is not a square in $L$. We claim that $\alpha = \alpha_0$ as above does the job. Indeed, suppose there is a $\gamma \in L$ such that $\gamma^2 = \alpha_0$. Then $N_{L/K}(\gamma)^2 = N_{L/K}(\alpha_0)$, contradicting our choice of $\alpha_0$.
\end{proof}

This classifies the structure of lattices over such fields to a satisfactory extent. The following is also an immediate corollary of our analysis.

\begin{cor}\label{ct7.1}
Let $K$ be a field with a linear lattice over it. Then either $K$ is real closed or any finite division algebra over it is a field extension. 
\end{cor}

\section{Quaternion Algebras and Merkurjev's Theorem} \label{s6}

In this section, we state some results about Brauer groups and quaternion algebras which play a crucial role in our proof of the targeted result.
\begin{defn}\label{d:qalg}
Let $K$ be a field of characteristic unequal to $2$. We say $D/K$ is a \textit{Quaternion Algebra} if there are $a,b\in K$ and $i, j, k$ in $D$ such that
\begin{enumerate}
    \item $D$ has a $K-$basis $1,i,j,k$.
    \item $i^2=a, j^2 = b, ij = -k, ji = k$.
\end{enumerate}

For a quaternion algebra, we have the \textit{norm map} $N(t+xi+yj+zk)=t^2 - ax^2 - by^2+abz^2$. Denote the above quaternion algebra as $Q(a,b)$.
\end{defn}

We state a few properties of such algebras. For a proof, one can see \cite{voi21}.

\begin{prop}\label{p654} 
let $K$ be a field with characteristic unequal to $2$.
\begin{enumerate}
\item Every quaternion algebra over $K$ is isomorphic to either a division algebra with centre $K$ or is isomorphic to $M_2(K)$.
\item The former alternative holds if and only if the norm map is isotropic as a quadratic form over $K$.
\item Every degree $2$ $K-$central division algebra is a quaternion algebra. 
\end{enumerate}
\end{prop}

As a consequence, we derive the following proposition. 

\begin{prop}\label{p8}
    Let $K$ be a perfect field such that
    \begin{enumerate}
        \item $(-1)$ is a square in $K$.
        \item $G = G(/K)$ is pro-2 and non-trivial.
        \item For every non-trivial finite extension $L/K$, $Br(L) = 0$.
    \end{enumerate}
    Then $Br(K) = 0$.
\end{prop}

\begin{proof}
	If characteristic is $2$, then $K=K^2$. Let $L/K$ be a cyclic extension of degree $2$. Then by Corollary \ref{cp2.2}, $Br(L)=0$ would imply that $Br(K)$ is $0$. So we may assume characteristic of $K$ is not $2$.

    Suppose $Br(K) \neq 0$. Let $L/K$ be a minimal proper extension. If $L=K'$, then we would immediately lead to a contradiction by Theorem \ref{t6}. So let us assume that there is another minimal proper extension $L_1/K$.

    Now, let $D$ be a non-commutative finite division algebra with center $K$. This will exist as $Br(K)$ is assumed to be nonzero. As $L/K$ is a splitting extension of $D$, we get that $D$ has degree $2$ over $K$. So $D$ is a quaternion algebra, say $D\cong Q(a,b)$ with $a,b\in K^*$. Let $K_1  = K(\sqrt{a})$ and $K_2=K(\sqrt{b})$. If $K_1=K_2$, then $a=m^2b$ for some $m\in K^*$. In that case, the norm function for $Q(a,b)$ becomes $N(t+xi+yj+zk)=t^2 - m^2bx^2-by^2+m^2b^2z^2$. As $(-1)$ is a square, we realise that $(\sqrt{-1}mb + k)$ is a non-zero element with a zero norm, hence $Q(a,b)$ is not a division algebra. 

    So we must have that $K_1,K_2$ are linearly disjoint over $K$. Thus $N_{K_1K_2/K_1}$ is surjective as $Br(K_1)=0$. In other words, there are $t,x,y,z\in K$ such that 
    $\sqrt{a} = (t+x\sqrt{a})^2 - b (y+z\sqrt{a})^2$. By using the linear independence of $1$ and $\sqrt{a}$ over $K$, we get that $t^2+x^2a - by^2 - abz^2 = 0$ and $1 = 2tx-2byz$. Hence, $\alpha = t + \sqrt{-1}xi +bj+\sqrt{-1}z$ is a non-zero element with a zero norm, showing once again that $Q(a,b)$ is not a division algebra.

   This contradicts that $Br(K)\neq 0$. Hence $Br(K)=0$, as desired. 
\end{proof}

In \cite{bru68}, Brumer and Rosen made the following conjecture.

\textit{\textbf{Conjecture.} 
    Let $K$ be a field and $p$ be a prime unequal to the characteristic of $K$. Then $B_p$, the $p-$part of $B=Br(L/K)$ satisfies one of the following.
    \begin{enumerate}
        \item It is divisible (possibly trivial).
        \item $p=2$ and it is an elementary abelian $2-$group. 
    \end{enumerate}
}

Merkurjev proves the following theorem in \cite{mer83}. 

\begin{thm}[Merkurjev]\label{t8} 
Suppose $K$ is a field and $p\neq char(K)$ be a prime. Let $\mu_p$ be the set of all $p^{th}$ roots of unity in some algebraic closure of $K$. If $[K(\mu_p):K]\leq 3$, then Brumer-Rosen conjecture holds. 
\end{thm}

We use this theorem to derive the following crucial result.

\begin{prop}\label{p9}
    Let $K$ be a perfect field for which $Br(K)$ is a finite non-trivial group and $Br(L) = 0$ for every nontrivial finite extension $L/K$. Then $G(/K)$ is a pro-$2$ group and $(-1)$ is not a square in $K$.
\end{prop}

\begin{proof}
	Let $B=Br(K)$, $p$ be a prime and $B_p$ denote the Sylow $p-$part of $B$. Suppose $B_p\neq 0$. Let $L = K(\mu_p)$. As $[L:K]\leq p-1$, Proposition \ref{p1}, parts (ii) and (iv) would yield an injection $B_p \hookrightarrow Br(L)$. If $L\neq K$, then $Br(L)=0$ implies $B_p=0$. Now suppose $L=K$. By Merkurjev's theorem, if $p$ is an odd prime, then $B_p$ is divisible and finite. Hence, $B_p=0$.  
	
	Therefore, $B$ is an elementary abelian $2-$group. Now let $L/K$ be a proper extension. Clearly, $Br(K)=Br(L/K)$ as $Br(L)=0$. So $Br(K)$ is $[L:K]$ torsion and also $2-$torsion. As $B[K]\neq 0$, this implies that $2$ divides $[L:K]$. So, by the proof of Proposition \ref{p4:ac}, we get that $G(/K)$ is pro-$2$. Also, Proposition \ref{p8} would show that $(-1)$ is not a square in $K$.
\end{proof}

The next theorem studies the only remaining case, that is, when $G(/K)$ is pro-2, $Br(K)$ is non-zero and $-1$ is not a square in $K$

\begin{thm}\label{t9} 
Let $K$ be a perfect field such that $(-1)$ is not a square in $K$, $G=G(/K)$ is a pro-2 group and $Br(K)\neq 0$. Then $K$ is either real closed or has a quadratic extension $M/K$ such that $Br(M)\neq 0$.
\end{thm}

\begin{proof}
Throughout this proof, we fix an algebraic closure of $K$ and denote by $i$ a square root of $(-1)$ in it. Suppose $K$ is not real closed. If $L=K(i)$ has a nonzero Brauer group, we are done. Otherwise, we must have that $Br(L/K)=Br(K)\neq 0$. Hence, the norm map $N_{K(i)/K}$ is not surjective 

Suppose there is a $T$ in $K$ such that $T$ is not in the image of the norm and $K(\sqrt{T})\neq K(i)$. Then we claim that $\sqrt{-T}$ is not in the image of $N_{M(i)/M}$ where $M= K(\sqrt{-T})$. Indeed, if it was, then there will be $a,b,c,d\in K$ such that $(a+b\sqrt{-T})^2+(c+d\sqrt{-T})^2 = \sqrt{-T}$. That would imply, by the linear independence of $1,\sqrt{-T}$ over $K$, that $a^2+c^2 = (b^2+d^2)T$ and $2(ab+cd) = 1$. If $b^2+d^2 = 0$, then $a^2+c^2 = 0$. For $s\in K$, we get $s = (a+bs)^2+(c+ds)^2$, showing that $N_{K(i)/K}$ is surjective. This is a contradiction.

Otherwise, we will get $T = \frac{(a^2+c^2)(b^2+d^2)}{(b^2+d^2)^2} = \frac{(ab-cd)^2}{(b^2+d^2)^2} + \frac{(ad+bc)^2}{(b^2+d^2)^2}$, whence we get a contradiction to the choice of $T$. So, $\sqrt{-T}$ is not in the image of the norm $N_{M(i)/M}$. This would imply that $Br(M(i)/M)\neq 0$, which, in turn, would show that $Br(M)\neq 0$.

Finally, if there is no such $T$, then for each $u$ not in the image of the norm $N_{L/K}$, we get that $u=-a^2$ for some $a\in K$. So every element of $K$ is either a square or the negative of a square. This would show that $L=K'$. As $N_{L/K}$ is not surjective, Proposition \ref{p7} would imply that $K$ is real closed. 
\end{proof}

\section{The Main Result} \label{s7}

Finally, we are in a position to state and prove the main result of this article.

\begin{thm}\label{t10} 
A field $K$ having only finitely and positively many nonisomorphic noncommutative finite division algebras over it is real closed.
\end{thm}

\begin{proof}
By Corollary \ref{p1.1}, we realise that if $K$ has finitely many noncommutative division algebras over it, then either $Br(L)=0$ for each $L/K$ finite or there is a maximal finite extension $M/K$ such that $Br(M)$ is non-zero, finite and for each $L/M$ non-trivial, $Br(L)=0$. In the latter case, by Proposition \ref{p9}, $G=G(/M)$ is pro-$2$ and $(-1)$ is not a square root in $M$. By Theorem \ref{t9}, $M$ is either real closed or it has a quadratic extension with non-trivial Brauer group. The latter alternative contradicts the choice of $M$. Thus, $M$ must be real closed. But a real closed field has no finite sub-extension by the Artin-Schreier Theorem. This shows that $K=M$, settling our claim.
\end{proof}

In other words, for any field $K$, the cardinality of $S(K)$ is either $0$, $1$ or infinite, with $1$ occurring if and only if $K$ is real closed.

\section*{Acknowledgements}

I would like to thank Professor James Borger of the Australian National University for initiating the discussions leading to this article, reviewing my work and also suggesting changes, approaches and further generalisations. I would like to thank the Australian National University and the organisers of the Future Research Talent (FRT) program as most the work presented was done during my stay there. This work has been supported partly by the FRT award. I would like to thank Professor L. H. Rowen of the Bar-Ilan University for his suggestions about improving the structure of this article. I would also like to thank my professors and colleagues at my home institute, Indian Statistical Institute, Kolkata for their guidance, support and for nominating me for the aforementioned program.

\end{document}